\documentclass[a4paper, 12pt]{amsart}
\usepackage[utf8]{inputenc}

\usepackage[el,nf]{coelacanth}
\usepackage[T1]{fontenc}

\let\oldnormalfont\normalfont
\def\normalfont{\oldnormalfont\mdseries}
\usepackage{pgfplots}
\pgfplotsset{compat=1.15}
\usepackage{mathrsfs}
\usetikzlibrary{arrows}

\usepackage{xfrac}

\usepackage{color}
\usepackage[colorlinks=true,linkcolor=colorref,citecolor=colorcita,urlcolor=colorweb]{hyperref}
\definecolor{colorcita}{RGB}{21,86,130}
\definecolor{colorref}{RGB}{5,10,177}
\definecolor{colorweb}{RGB}{177,6,38}

\usepackage{amsmath}
\usepackage{amsthm}
\usepackage{xcolor}
\usepackage{soul}

\newtheorem{thm}{Theorem}[section]

\newtheorem{defi}[thm]{Definition}

\DeclareFontFamily{U}{mathx}{\hyphenchar\font45}
\DeclareFontShape{U}{mathx}{m}{n}{
      <5> <6> <7> <8> <9> <10>
      <10.95> <12> <14.4> <17.28> <20.74> <24.88>
      mathx10
      }{}
\DeclareSymbolFont{mathx}{U}{mathx}{m}{n}

\DeclareMathAccent{\widecheck}{0}{mathx}{"71}
\newcommand{\zC}{\mathbb{C}}
\newcommand{\zN}{\mathbb{N}}
\newcommand{\zR}{\mathbb{R}}

\newcommand{\ve}{\varepsilon}

\newcommand{\cc}[3]{\mathbf{c}_\mathfrak{#3}({#1}, {#2})}
\newcommand{\ccl}[2]{\mathbf{c}_\mathfrak{#2}({#1})}

\newcommand{\vertiii}[1]{{\left\vert\kern-0.25ex\left\vert\kern-0.25ex\left\vert #1 
    \right\vert\kern-0.25ex\right\vert\kern-0.25ex\right\vert}}

\begin{document}

\title[A Banach space with constant of analyticity less than one]{A $2$-dimensional real Banach space with constant of analyticity less than one}

\author[Jorge Tom\'as Rodr\'{i}guez]{Jorge Tom\'as Rodr\'{i}guez}
\address{Departamento de Matem\'{a}tica and NUCOMPA, Facultad de Cs. Exactas, Universidad Nacional del Centro de la Provincia de Buenos Aires, (7000) Tandil, Argentina and CONICET\newline
	\href{https://orcid.org/0000-0003-4693-2498}{ORCID: \texttt{0000-0003-4693-2498} } }
\email{jtrodriguez@nucompa.exa.unicen.edu.ar}

\begin{abstract} We show that on the real $2$-dimensional Banach space $\ell_1^2$ there is an analytic function $f:B_{\ell_1^2}\rightarrow \zR$ such that its power series at origin has radius of uniform convergence one, but for some $a\in B_{\ell_1^2}$ the power series centred at that point has radius of uniform convergence strictly less than $1-\|a\|$. This result highlights a fundamental distinction in real analytic functions (compared to complex analytic functions), where the radius of analyticity can differ from the radius of uniform convergence. Moreover, this example provides the first non-trivial upper bound for the constant of analyticity.
\end{abstract}

\thanks{This work was partially supported by 
CONICET PIP 11220200101609CO,   %vero y tomas
and ANPCyT PICT 2018-04250. %tomas
}

\maketitle

\section{Introduction}

In this work we study the radius of analyticity of real analytic function, that is, of those functions that can be approximated by polynomials. To establish the foundation for our discussion, we will first introduce some essential concepts briefly.

A function $P:X\to E$ between Banach spaces is a \textbf{continuous $k$-homogeneous polynomial} if there exists a (unique) continuous symmetric $k$-linear mapping $\widecheck{P}$ such that \break $P(x) = \widecheck{P}(x,\dots, x)$. We write $\mathcal{P}(^kX,E)$ for spaces of  $k-$homogeneous polynomials. When $E=\mathbb K$, the scalar field considered, we omit it from the notation and write  $\mathcal{P}(^kX)$. The norm of a polynomial is defined as follows
$$\|P\|=\sup\{ \|P(x)\|_E: \|x\|_X\leq 1\}.$$
It's important to note that our main result in this work pertains to a finite-dimensional space. Thus, it is worth mentioning that in this setting, the above definition of $k-$homogeneous polynomial coincides with the usual definition on finite variables. 

Given an open subset $U\subseteq X$, a function $f:U\rightarrow E$ is \textbf{analytic} if for every $a\in U$ there is $r>0$ and a sequence of polynomials $P_k \in \mathcal{P}(^kX,E)$ (depending on $f$ and $a$) such that
$$f(x)=\sum_{k=1}^\infty P_k(x-a)$$
uniformly on $B(a,r)$, the ball of centre $a$ and radius $r$. For a more detailed exposition on real analytic functions, we refer to \cite{hajek2014smooth}, while for further information on complex analytic functions (also known as holomorphic functions), we refer to \cite{chae2020holomorphy, Din99, Muj86}.

The series
$$S(f,a)=\sum_{k=1}^\infty P_k(x-a)$$
is called the  \textbf{power series} associated to $f$ about $a$. The \textbf{radius of uniform convergence} is defined as
$$ R(f,a)=\sup\{ r : S(f,a) \text{ converges uniformly in } B(a,r)\}.$$
It is well known that 
$$R(f,a)=\liminf_{k\to \infty} \frac{1}{\|P_k\|^{\frac{1}{k}}}.$$
We always assume that $B(a, R(f,a))$ is included in $U$, the domain of $f$. If not, we can extend $f$  to $U\cup B(a, R(f,a))$ and work with this extension.

Next define the \textbf{radius of analyticity} of $f$ about $a$
\begin{eqnarray*}
    R_A(f,a)&=&\sup\{ r : S(f,b)\text{ converges uniformly in } B(b,r-\|b-a\|) \forall b \in B(a,r)\}\\
    &=&\inf\{ \|b-a\|+R(f,b): b\in B(a,R(f,a))\} \
\end{eqnarray*}
In a typical introductory complex analysis course, it is established that $R(f,a)=R_A(f,a)$, rendering the latter definition irrelevant within this specific context. This outcome can be generalized to encompass arbitrary complex Banach spaces, thanks to the application of Cauchy's estimates. This conclusion also holds true when dealing with real Hilbert spaces, using that the polarization constant of these spaces is one. However, for general real Banach spaces, we do not have Cauchy estimates, and the polarization constant is one only for Hilbert spaces. This leads to a natural question: Can we find an example for which \break $R(f,a)>R_A(f,a)$? This question motivates the next definition (see \cite{boyd2018radius, nguyen2009lower}).

\begin{defi}
    Let $X$ be a real Banach space, the \textbf{constant of analyticity} $\mathcal{A}(X)$ of $X$ is the supremum of the positive real numbers $\rho$ such that every power series at the origin in $X$, with unit radius of uniform convergence, has radius of analyticity $\rho$. That is
    $$\mathcal{A}(X) = \inf\{R_A(f,0)\},$$
    where the infimum is taken over all the analytic series given by a power series about the origin with unit radius of uniform convergence.
\end{defi}

Whether or not we always have $R(f,a)=R_A(f,a)$ is a long-standing question (see \cite[Problem 1.177]{hajek2014smooth}). This is equivalent to ask whether $\mathcal{A}(X)$ is always one or not. In \cite{taylor1938additions}, A. E. Taylor showed that
$$\frac{1}{e\sqrt{2}}\leq \mathcal{A}(X).$$
This lower bound was improved later. The current best general lower bound is due to P. Hajek and M. Johanis
$$\frac{1}{\sqrt{2}}\leq \mathcal{A}(X).$$
Independently in \cite{papadiamantis2016radius}, M. K. Papadiamantis and Y. Sarantopoulos proved the same lower bound with different techniques.

In this work, we show, using complexification techniques, that for some spaces $\mathcal{A}(X)<1$. In \cite{boyd2018radius} the authors show that for any real Banach space 
$$\mathcal{A}(\ell_1(\zR)) \leq \mathcal{A}(X).$$
Thus, the natural candidate to prove our result is $\ell_1(\zR)$, although the $2$-dimensional space $\ell_1^2(\zR)=(\zR^2, \|\cdot\|_1)$ will be sufficient.

For further differences between real and complex Banach spaces, we refer the reader to the survey \cite{moslehian2022similarities} and the references therein.

%%%%%%%%%%%%%%%%%%%%%%%%%%%%%%%%%%%%%%%%%%%%%%%%%%%%%%%%%%%%%%%%%%%%%%%%%%
%%%%%%%%%%%%%%%%%%%%%%%%%%%%%%%%%%%%%%%%%%%%%%%%%%%%%%%%%%%%%%%%%%%%%%%%%%
%%%%%%%%%%%%%%%%%%%%%%%%%%%%%%%%%%%%%%%%%%%%%%%%%%%%%%%%%%%%%%%%%%%%%%%%%%

%%%%%%%%%%%%%%%%%%%%%%%%%%%%%%%%%%%%%%%%%%%%%%%%%%%%%%%%%%%%%%%%%%%%%%%%%%
%%%%%%%%%%%%%%%%%%%%%%%%%%%	%%%%%%%%%%%%%%%%%%%%%%%%%%%%%%%%%%%%%%%%%%%%%%%
%%%%%%%%%%%%%%%%%%%%%%%%%%%%%%%%%%%%%%%%%%%%%%%%%%%%%%%%%%%%%%%%%%%%%%%%%%

\section{Main result}\label{section main}
Let us consider the space $\ell_1^2(\zR)$. The space $\ell_1^2(\zC)$  is a complexification of $\ell_1^2(\zR)$. Moreover, it is the complexification obtained using Bochnak's procedure. To see the definition of complexification procedures, we refer the reader to \cite[Section 2]{munoz1999complexifications}. Given $P\in \mathcal{P}(^k\ell_1^2(\zR))$ there is one extension $\Tilde{P}\in \mathcal{P}(^k\ell_1^2(\zC))$ (see \cite[Section 3]{munoz1999complexifications}). For any natural number $k\in \zN$ we define
\begin{eqnarray*}
\cc{k}{\ell_1^2(\zR)}{b}&=&\inf\{M >0: \Vert \widetilde{P} \Vert \leq M \Vert P\Vert, \text{ for all } P\in\mathcal{P}(^k\ell_1^2(\zR))\}\\
&=& \sup\{ \Vert \widetilde{P} \Vert : P\in\mathcal{P}(^k\ell_1^2(\zR)) \text{ with } \|P\|\leq 1\}.\
\end{eqnarray*}
We call $\cc{k}{\ell_1^2(\zR)}{b}$ the $k^{th}$ \textbf{complexification constant} of the Banach space $\ell_1^2(\zR)$ with \break Bochnak's complexification procedure. We also define
$$\ccl{\ell_1^2(\zR)}{b}=\limsup_{k\to\ \infty} (\cc{k}{\ell_1^2(\zR)}{b})^{\frac{1}{k}},$$
the \textbf{complexification constant}. For a more thorough introduction to complexification constant and their basic properties, we refer the reader to \cite{rodinpress}.

Although we do not know the exact value of $\ccl{\ell_1^2(\zR)}{b}$, in \cite[Proposition 2.6]{dimant2022polarization} the following estimation was obtained
$$\sqrt[4]{2} \leq \ccl{\ell_1^2(\zR)}{b}\leq \sqrt{2}.$$
Using the complexification constants as a tool, we can state and prove our main result.
\begin{thm} For the real Banach space $\ell_1^2(\zR) $ we have
    $$\frac{1}{\ccl{\ell_1^2(\zR)}{b}} \leq\mathcal{A}(\ell_1^2(\zR))\leq \frac{1}{2}+\frac{1}{2\ccl{\ell_1^2(\zR)}{b}} <1.$$
\end{thm}

\begin{proof}
    The lower bound can be established through a straightforward complexification argument, taking advantage of the fact that on complex Banach spaces, the radius of uniform convergence and the radius of analyticity coincide.

    Let us prove the upper bound. Given $\ve>0$ we will construct a power series at the origin with unit radius of uniform convergence but with radius of analyticity less than 
    $$\left(\frac{1}{2}+\frac{1}{2\ccl{\ell_1^2(\zR)}{b}} \right) \frac{1}{1-\ve}.$$
    
First fix a norm one polynomial $P\in \mathcal{P}(^k\ell_1^2(\zR))$ such that
$$\sqrt[k]{\|\tilde{P}\|} > \ccl{\ell_1^2(\zR)}{b} (1-\ve).$$
Next take a norm one vector $(\alpha, \beta)\in \ell_1^2(\zC)$ for which
\[ |\tilde{P}(\alpha, \beta)| =\|\tilde{P}\|.\]
We may assume $|\alpha| \geq \frac{1}{2}\geq |\beta|$. Also, replacing $(\alpha, \beta)$ by $e^{i\theta}(\alpha, \beta)$ for some suitable $\theta$ we may suppose $\alpha \in \zR_{>0}$. 

We additionally assume that the argument of the complex number $\tilde{P}(\alpha, \beta)$ is rational (modulus $2\pi$), which implies the existence of an integer $m\in \zN$ such that $\tilde{P}(\alpha, \beta)^m\in \zR_{>0}$. This simplifies the proof considerably. At the end of the proof, we will provide guidance on how to handle cases where this assumption does not hold.

Consider the unit ball $B(0,1)\subset \ell_1^2(\zR)$ and the analytic function $f:B(0,1)\rightarrow \zR$ given by 
$$f(x)=\sum_{n=1}^\infty P^{mn}(x).$$
We have
$$R(f,0)= \liminf_{n \to \infty}\frac{1}{\sqrt[mnk]{\| P^{mn}\|}}= \frac{1}{\sqrt[k]{\| P\|}}=1.$$

Take $\alpha'=\frac{\alpha}{\sqrt[k]{\|\tilde{P}\|}}< \frac{\alpha}{\ccl{\ell_1^2(\zR)}{b} (1-\ve)}$, and consider the power series of $f$ at $(\alpha',0)$
$$f(x)= \sum_{j=1}^\infty Q_j(x-(\alpha',0)).$$
It is sufficient to show that
$$R(f,(\alpha',0)) \leq \frac{|\beta|}{(1-\ve)},$$
since this implies
\begin{eqnarray*}
\mathcal{A}(\ell_1^2(\zR)) &\leq& R_A(f,0)\\ 
 &\leq&  \| (\alpha',0)- 0\|+R(f, (\alpha',0))\\
&\leq&  \alpha'+\frac{|\beta|}{(1-\ve)}\\
&<&\frac{\alpha}{\ccl{\ell_1^2(\zR)}{b} (1-\ve)}+\frac{|\beta|}{(1-\ve)}\\
&\leq& \left(\frac{1}{2\ccl{\ell_1^2(\zR)}{b} }+\frac{1}{2}\right)\frac{1}{1-\ve},\
\end{eqnarray*} 
where in the last inequality, we use that $|\alpha| \geq \frac{1}{2}\geq |\beta|$ and that $|\alpha|+|\beta|=1$.

Now let us consider the complexification of $f$
$$\tilde{f}(z)= \sum_{n=1}^\infty \widetilde{P}^{mn}(z).$$
Then
$$R(\tilde{f},0)= \liminf_{n \to \infty}\frac{1}{\sqrt[mnk]{\| \tilde{P}^{mn}\|_\mathfrak{b}}}= \frac{1}{\sqrt[k]{\| \tilde{P}\|_\mathfrak{b}}}.$$
In particular $\tilde{f}$ is defined on $B\left(0, \frac{1}{\sqrt[k]{\| \tilde{P}\|_\mathfrak{b}}}\right)$. If we consider the power series centred at $(\alpha',0)$, we have 
$$\tilde{f}(z)= \sum_{j=1}^\infty \widetilde{Q_j}(z-(\alpha',0)),$$ and
\begin{eqnarray}
R(f,(\alpha',0)) &=&\liminf_{j \to \infty}\frac{1}{\sqrt[j]{\| Q_j\|}}\nonumber \\ 
& \leq &\liminf_{j \to \infty}\frac{\sqrt[j]{\cc{j}{\ell_1^2(\zR)}{b}}}{\sqrt[j]{\|\widetilde{ Q_j}\|_\mathfrak{b}}} \label{eq Q} \\ 
& \leq &\limsup_{j \to \infty}\sqrt[j]{\cc{j}{\ell_1^2(\zR)}{b}}\,\liminf_{j \to \infty}\frac{1}{\sqrt[j]{\|\widetilde{ Q_j}\|_\mathfrak{b}}} \nonumber \\ 
&=&\ccl{\ell_1^2(\zR)}{b}R(\tilde{f},(\alpha',0)).\nonumber\
\end{eqnarray}
Thus, if we obtain an upper bound for $R(\tilde{f},(\alpha',0))$ we get one for $R(f,(\alpha',0))$.

Let us define $\beta'=\frac{\beta}{\sqrt[k]{\|\tilde{P}\|}}< \frac{\beta}{\ccl{\ell_1^2(\zR)}{b} (1-\ve)}$. We claim that $\tilde{f}$ can not be extended continusly to $(\alpha',\beta')$. Indeed, if $t\in (0,1)$ we have
\begin{eqnarray*}
    \tilde{f}(t(\alpha',\beta'))&=&\sum_{n=1}^\infty t^{kmn} \widetilde{P}^{mn}(\alpha',\beta') \\
    &=&\sum_{n=1}^\infty t^{kmn}\\
    &=&\frac{1}{1-t^{km}}-1\xrightarrow[t\to 1]{} +\infty.\
\end{eqnarray*}

Then, since $\tilde{f}$ can not be extended to $(\alpha',\beta')$, we have
$$R(\tilde{f},(\alpha',0))\leq \|(\alpha',0)-(\alpha',\beta')\|=|\beta'|< \frac{\beta}{\ccl{\ell_1^2(\zR)}{b} (1-\ve)}.$$
This gives the desired upper bound
\begin{eqnarray*}
R(f,(\alpha',0)) & \leq &\ccl{\ell_1^2(\zR)}{b}R(\tilde{f},(\alpha',0))  \\ 
&\leq &\ccl{\ell_1^2(\zR)}{b} \frac{\beta}{\ccl{\ell_1^2(\zR)}{b} (1-\ve)}  \\
&= & \frac{|\beta|}{(1-\ve)}.\
\end{eqnarray*}

We only need to deal with the case $\tilde{P}(\alpha, \beta)$ has irrational argument. If that is the case, right after we fix $\alpha$ and $\beta$ we take $\theta>0$ such that $\tilde{P}(e^\theta(\alpha, \beta))$ has a rational argument and  small enough such that 
$$\|(\alpha',0)-(e^\theta\alpha',e^\theta\beta')\|=|\alpha'|(1-e^\theta)|+|\beta'|< \frac{\beta}{\ccl{\ell_1^2(\zR)}{b} (1-\ve)}.$$
We can do this because $|\beta'|$ is strictly smaller than $\frac{\beta}{\ccl{\ell_1^2(\zR)}{b} (1-\ve)}$.
Then we take $m$ such that $\tilde{P}^m(e^\theta(\alpha, \beta))\in \zR_{>0},$ and proceed as before. Following the same steps as before, we  conclude that $\tilde{f}$ can not be extended to $(e^\theta\alpha',e^\theta\beta')$, giving
$$R(\tilde{f},(\alpha',0))\leq \|(\alpha',0)-(e^\theta\alpha',e^\theta\beta')\|< \frac{\beta}{\ccl{\ell_1^2(\zR)}{b} (1-\ve)},$$
which is exactly what we need to conclude the desired result.

\end{proof}

It's worth noting that in the preceding proof, if we possessed the explicit expression of $P$, this technique could potentially yield a more favourable upper bound. However, due to the absence of such information, we lack any insights into the polynomials $Q_j$. Therefore, on \eqref{eq Q} we utilize the bound
$$ \|\tilde{Q}_j \|\leq \cc{j}{\ell_1^2(\zR)}{b} \|Q_j\|. $$
But this bound may not be optimal. For instance, if we were to have $\|\tilde{Q}_j \|=\|Q_j\|$, the same proof would lead to $\ell_1^2(\zR)$ having constant of analyticity equal to $\frac{1}{\ccl{\ell_1^2(\zR)}{b}}$.

%\subsection*{Acknowledgments}

%%%%%%%%%%%%%%%%%%%%%%%%%%%%%%%%%%%%%%%%%%%%%%%%%%%%%%%%%%%%%%%%%%%%%%%%%%
%%%%%%%%%%%%%%%%%%%%%%%%%%%%%%%%%%%%%%%%%%%%%%%%%%%%%%%%%%%%%%%%%%%%%%%%%%
%%%%%%%%%%%%%%%%%%%%%%%%%%%%%%%%%%%%%%%%%%%%%%%%%%%%%%%%%%%%%%%%%%%%%%%%%%

%%%%%%%%%%%%%%%%%%%%%%%%%%%%%%%%%%%%%%%%%%%%%%%%%%%%%%%%%%%%%%%%%%%%%%%%%%
%%%%%%%%%%%%%%%%%%%%%%%%%%%	%%%%%%%%%%%%%%%%%%%%%%%%%%%%%%%%%%%%%%%%%%%%%%%
%%%%%%%%%%%%%%%%%%%%%%%%%%%%%%%%%%%%%%%%%%%%%%%%%%%%%%%%%%%%%%%%%%%%%%%%%%

\bibliography{biblio.bib}

\bibliographystyle{abbrv}

\end{document}